# A Voronoi Diagram-Based Approach for AC Optimal Power Flow

Mohammed N. Khamees, Kai Sun

Electrical Engineer and Computer Science Department, University of Tennessee, Knoxville, TN, USA
mkhamees@vols.utk.edu, kaisun@utk.edu

**Abstract:** The primary goal of Optimal Power Flow (OPF) is to optimize the operation of a power system while meeting the demand and adhering to operational constraints. This paper presents a new approach for AC OPF. First, the approach constructs a Voronoi diagram by distributing multiple sample points representing potential solutions throughout the search space. Then, it recursively adds new sample points including a tentative optimal point from the continuous gradient-projection method, a point in the most sparsely populated region to ensure high fidelity, and the connecting point, until the stopping criterion is met. The proposed approach is illustrated in detail using the IEEE 9-bus system and then validated on the IEEE 39-bus and 118-bus systems to verify the quality of the obtained solution.

## 1. Introduction

### 1.1. Background

Optimal power flow (OPF) is an optimization problem that determines the most efficient operational state of a power system, considering a specified objective function while adhering to operational constraints and meeting demand requirements. The inherent nonlinear relationships between power system variables and parameters, such as complex powers and voltage magnitudes as well as the constraints and the objective function of the OPF problem, contribute to its nonconvex nature [1], which significantly increases the likelihood of encountering local OPF solutions, making the reliable computation of the global optimal solution a challenging task [2], [3].

### 1.2. Literature review

Since Carpentier's seminal work in 1962 [4], various techniques have been employed to address OPF problems, including Gradient methods, Newton's method, sequential programming, and interior point methods. A comprehensive review of these techniques is provided in [5]. However, these methods can suffer from convergence failures or being trapped by local optima, making global optimality uncertain [2], [5], [6]. Over the past two decades, convex relaxation techniques, such as semidefinite programming (SDP) [1] and second-order cone programming (SOCP) [7] have gained attention. While these methods enable accurate computations for large-scale networks with global guarantees if the relaxation is exact, they provide a lower bound on the unknown globally optimal in case the relaxation is inexact. In the worst case they may fail to yield physically meaningful solutions [8]. For instance, most experiments conducted on a diverse dataset in [9] revealed a duality gap greater than zero, indicating inexact relaxation. Another limitation of the current SDP solvers lies in their precision. Even when the duality gap is zero, they do not guarantee the exact global optimum mathematically but only provide an approximation [10]. Using polynomial optimization theory, the moment/sum-of-squares relaxations can construct a stronger convex relaxation approach for addressing OPF problems. However, the size of the SDPs grows drastically with the number of OPF variables and with the order of the relaxation [11]. To reduce the computational burden of SDP relaxations, studies such as [11] and [12] exploited sparsity techniques to achieve 1% global optimality gaps. To address the limitation of SDP in ensuring a rank-one solution, [13] introduces an iterative approach that solves a sequence of semidefinite programs. The solutions progressively converge to a rank-one matrix, corresponding to a stationary point. However, the resulting solutions still exhibit small optimality gaps.

Ref [14] employs a gradient-based approach to explore the relationship between KKT saddle points and optimal solutions. Identifying all optimal solutions requires computing all KKT saddle points, which can result in thousands of KKT points, even for a problem of moderate size. The branch tracing method and monotone search strategy have been proposed to identify multiple local OPF solutions [15]. These methods compute stationary solutions using the first-order optimality conditions and apply the second-order sufficient conditions to classify local extrema. However, they require intricate formulations and are computationally intensive [15]. Metaheuristic approaches, including evolutionary algorithms and particle swarm optimization, offer the potential to theoretically converge to a global optimum. However, these methods typically demand substantial computational resources, often relying heavily on meticulous parameter tuning for effectiveness [16].

An alternative approach to achieving global optimality is the enumeration of local solutions. Ref. [2] identifies multiple local OPF solutions by employing local optimization algorithms initialized with numerous random starting points. This method demonstrates its capability to uncover several local optima across various IEEE system models where local algorithms frequently succeed in identifying the global optimum. Furthermore, it has been observed that SDP often fails in most cases where local optima are present [2].

### 1.3. Contributions

In this paper, a new approach is proposed to globally solve the AC OPF problem by exploiting the geometrical structure of the OPF search space through a computational geometry technique, the Voronoi diagram [17]. Toward global optimization, the Voronoi diagram is applicable to characterize the search space of any dimension but requires sufficient sample points. In this paper, the sample points are efficiently added by using the continuous-projected gradient method with enforced constraint restoration [18]. Specifically, the method constructs a fictitious nonlinear dynamical system that models the OPF solution process. The asymptotically stable equilibria of this system correspond to local OPF solutions, whereas the other OPF critical stationary points—such as maxima and saddle points with the OPF formulated as a minimization problem—cannot be asymptotically stable. To achieve global optimality, the Voronoi diagram provides a global, geometrical characterization of the search space by partitioning it into cells and refining it iteratively with new sampling points. These points are selected based on the geometry of the feasible set, as revealed by the diagram, and heuristic insights derived from the nonlinear dynamical system.

*1.4. Organization*

The paper is organized as follows. Section II proposes the new approach to AC OPF, which provides its overview, introduces the continuous-projected gradient method with enforced constraint restoration and Voronoi diagram to OPF, and presents detailed algorithms. Sections III provides the OPF formulation to implement the proposed approach in this paper. Section IV first illustrates the approach on small examples and then validates it using IEEE test systems. Finally, section V concludes the paper.

## 2. Proposed Approach to AC Optimal Power Flow

*2.1. Overview of the Approach*

The proposed approach incorporates the continuous-projected gradient (CPG) method into a Voronoi diagram. Utilizing a tentative optimum as an initial point, the fictitious dynamical system constructed by the CPG method is integrated within the tentative optimum identified region. Characterized by the asymptotic stability and the constraint restoration feature of the dynamical system, an initial point can be feasible, infeasible or not power flow solution.

For a general nonlinear OPF problem:
$$\min_{\mathbf{u}} \; f(\mathbf{u},\mathbf{y})$$
$$s.t \; \mathbf{e}(\mathbf{u},\mathbf{y}) = \mathbf{0} \quad (1)$$
$$\mathbf{h}(\mathbf{u},\mathbf{y}) \leq \mathbf{0}$$

where $\mathbf{u}$ and $\mathbf{y}$ are respectively the control and state vectors, $f(\mathbf{u},\mathbf{y})$ is a differentiable cost function, and $\mathbf{e}(\mathbf{u},\mathbf{y})$ and $\mathbf{h}(\mathbf{u},\mathbf{y})$ are the equality and inequality constraints, respectively, both of which are differentiable. This problem can be converted to an equality-constrained optimization problem by introducing new unconstrained slack variables in vector $\mathbf{s}$ as follows:
$$\mathbf{h}(\mathbf{u},\mathbf{y}) + \mathbf{s}^2 = \mathbf{0} \quad (2)$$

where $\mathbf{s}^2$ denotes the vector where each element in $\mathbf{s}$ is squared. The dimension of $\mathbf{s}$ is equal to $\mathbf{h}(\mathbf{u},\mathbf{y})$ dimension. Although this increases the number of critical points of the cost function relative to the constraints, these additional points cannot serve as minimal points [14]. Then, the OPF (1) can be expressed concisely as in (3), which shares the same optimal points of the original problem (1) [18].

$$\min_{\mathbf{x}} \; f(\mathbf{x})$$
$$s.t \; \mathbf{g}(\mathbf{x}) = \mathbf{0} \quad (3)$$

where $\mathbf{x} \epsilon \mathbb{R}^n$ is the vector of all variables and $f(\mathbf{x}): \mathbb{R}^n \to \mathbb{R}$ and $\mathbf{g}(\mathbf{x}): \mathbb{R}^n \to \mathbb{R}^m$. Assume the Linear Independent Constraint Qualification (LICQ) condition holds, which guarantees that the gradients of the constraints are linearly independent. This assumption is generally true for OPF problem [3], [15]. Then, the Jacobian matrix of $\mathbf{g}(\mathbf{x})$, denoted by $\mathbf{J_g(x)}$, has a full rank. Denote the cost function gradient by $\nabla f(\mathbf{x})$. The proposed approach takes the following steps:

**Step 1:** Construct an autonomous dynamical system in the form of $dx/dt=\psi(\mathbf{x})$ according to $\nabla f(\mathbf{x})$ and $\mathbf{J_g(x)}$.
**Step 2:** Select initial sample points and compute a penalized cost function denoted by $f_p(\mathbf{x})$.
**Step 3:** Establish the Voronoi diagram using all sample points.
**Step 4:** Identify the tentative optimum and the candidate Voronoi region.
**Step 5:** Identify and add new sample points.
**Step 6:** End the process if the maximum number of iterations is reached; otherwise, return to **Step 3**.

The following flow chart shows the proposed approach:

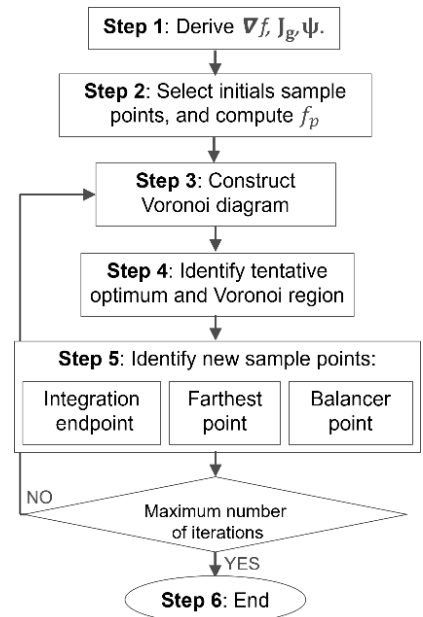

**Fig. 1**. Flow chart of the proposed approach.

The next subsections will present detailed algorithms of Step 1 to Step 5.

*2.2. Step 1: Construct an Autonomous Dynamical System*

In **Step 1**, the proposed approach uses the CPG method with enforced constraint restoration to construct a nonlinear dynamical system described by an ordinary differential equation characterizing the flows leading to local optima as introduced below [18].

$\mathbf{g}(\mathbf{x})$ in (3), let $m < n$, includes power flow equations and operational constraints. Its Jacobian has a rank of $m$ as calculated by:

$$\mathbf{J_g}(\mathbf{x}) = \left[\frac{\partial \mathbf{g}_i(\mathbf{x})}{\partial x_j}\right] \quad i \in \{1,..,m\}, j \in \{1,..,n\} \quad (4)$$

The OPF feasible space $\mathbf{FS_g}$ defined by (5) is an $(n-m)$-dimensional differentiable manifold in $\mathbb{R}^n$:

$$\mathbf{FS_g} \equiv \{\mathbf{x} \in R^n : \mathbf{g}(\mathbf{x}) = \mathbf{0}\} \quad (5)$$

Its tangent space ($\mathbf{T}_X$) at a feasible point $\mathbf{x} \in \mathbf{FS_g}$ can be defined by the null space of $\mathbf{J_g}(\mathbf{x})$:

$$\mathbf{N}(\mathbf{J_g}(\mathbf{x})) = \{\mathbf{z} \in R^n : \mathbf{J_g}(\mathbf{x})\mathbf{z} = \mathbf{0} \in R^m\} \quad (6)$$

The cost function gradient vector can be written as:

$$\nabla f(\mathbf{x}) = (\partial f / \partial x_1, ..., \partial f / \partial x_n)^T \quad (7)$$

Define the orthogonal projector $\mathbf{P}(\mathbf{x})$ into the tangent space as:

$$\mathbf{P}(\mathbf{x}) = \mathbf{I} - \mathbf{J_g^T}(\mathbf{x})(\mathbf{J_g}(\mathbf{x})\mathbf{J_g^T}(\mathbf{x}))^{-1}\mathbf{J_g}(\mathbf{x}) \quad (8)$$

which is an idempotent, symmetric matrix with its range perpendicular to its null [19]. Now, consider the following autonomous system for the OPF problem in (3):

$$\begin{bmatrix} \mathbf{I} & \mathbf{J_g^T}(\mathbf{x}) \\ \mathbf{J_g}(\mathbf{x}) & \mathbf{0} \end{bmatrix}\begin{bmatrix} \mathbf{\psi}(\mathbf{x}) \\ \mathbf{\lambda}(\mathbf{x}) \end{bmatrix} = -\begin{bmatrix} \nabla f(\mathbf{x}) \\ \mathbf{g}(\mathbf{x}) \end{bmatrix} \quad (9)$$

from which we have (10-11) with the $x$ argument is dropped:

$$\frac{dx}{dt} = \mathbf{\psi} = -\mathbf{J_g^T}(\mathbf{J_g}\mathbf{J_g^T})^{-1}\mathbf{g} - (I - \mathbf{J_g^T}(\mathbf{J_g}\mathbf{J_g^T})^{-1}\mathbf{J_g})\nabla f \quad (10)$$

$$\mathbf{\lambda} = (\mathbf{J_g}\mathbf{J_g^T})^{-1}(-\mathbf{J_g}\nabla f + \mathbf{g}) \quad (11)$$

$\mathbf{\psi}(\mathbf{x})$ will be explored to characterize the relation between its equilibria and the optimal solutions of problem (3).

There are the following two theorems:

**Theorem 1** [18]: A point $\mathbf{x}^* \in \mathbf{FS_g}$ is a critical point of the cost function $f(\mathbf{x})$ if and only if it is an equilibrium point of $\mathbf{\psi}(\mathbf{x})$, i.e., $\mathbf{\psi}(\mathbf{x}^*) = 0$.

**Theorem 2** [18]: If a critical point $\mathbf{x}^*$ is an isolated (regular) minimal point of (3) then it is asymptotically stable. Points which are not minimal points cannot be asymptotically stable points of $\mathbf{\psi}(\mathbf{x})$.

Theorem 1 establishes the connection between the critical points of OPF (3) and the equilibrium points of $\mathbf{\psi}(\mathbf{x})$. Theorem 2 asserts the relation between OPF and $\mathbf{\psi}(\mathbf{x})$, specifically in that each optimum of OPF is asymptotically stable equilibrium of $\mathbf{\psi}(\mathbf{x})$. Simultaneously, other critical points in OPF cannot serve as asymptotically stable points. Hence, $\mathbf{\psi}(\mathbf{x})$ defines a vector field in which the set of optima serves as its attractors. In ill-conditioned OPF problem where LICQ may fail, the CPG formulation can be adapted by employing alternative optimality conditions such as the Fritz–John conditions, ensuring robustness even when the KKT framework is not valid.

### 2.3. Step 2: Select Initial Sample Points and Computing a Penalized Cost Function

The penalized cost function $f_p(\mathbf{x})$ constructed in **Step 2** comprises the original cost function augmented by a penalty term proportional to the magnitude of constraint violations, which becomes zero when no violations are present. Solving an OPF problem involves identifying control variables, such as P-V bus voltage magnitude, transformer tap Settings and active power output, while state variables and functional constraints are determined subsequently based on the control variables. Therefore, for an OPF problem with $N$ control variables, the search space is $N$-dimensional. Next, using the control variables, which constitute the elements of a sample point, the power flow equations are solved to compute the state variables and functional constraints. Sample points can be classified into three categories based on feasibility:

I. **Violating equality constraints**: These points are not power flow solutions. In this case, the penalized cost function $f_p(\mathbf{x})$ is assigned to values larger than the points in the search space, to discourage selection of such points.

II. **Satisfying equality constraints but not inequality constraints:** These points are power flow solutions violating the inequality constraints. Their penalized cost function is computed incorporating the penalty proportional to the size of each violation $\alpha_k$ at a constant price of $c$:

$$f_p(\mathbf{x}) = f(\mathbf{x}) + c\sum_{k=1}^{K}\alpha_k \quad (12)$$

III. **Satisfying all constraints:** These points form feasible space (5) with no penalty applied.

### 2.4. Step 3: Establish a Voronoi Diagram

The Voronoi diagram is a computational geometry approach applicable to spaces of any dimension [17]. In this paper, the proposed approach will use it to divide the OPF search space into regions according to the distance from a set of sample points, which are composed of OPF control variables, e.g., bus voltages magnitude and power injections, offering geometric insights, such as identifying the most sparsely populated areas in space. It has been utilized to achieve global accuracy in function approximation for engineering design optimization [20], [21].

Given a set of $b$ distinct sampling points $\mathbb{P} = \{\mathbf{p}_1, ..., \mathbf{p}_b\}$ distributed across an $N$-dimensional space, a Voronoi diagram partitions the plane into $b$ cells, each enclosing the region closest to a specific point, as illustrated by Fig. 2 for a two-dimensional plane. There are two steps to construct a Voronoi diagram: 1) connect each sampling point to its neighboring points with lines (blue lines in the figure); 2) draw perpendicular bisectors between the connected pairs of points (black lines). This process generates a series of convex polygons (or convex polyhedra in a higher-dimensional space), which are called Voronoi regions. Each Voronoi region $V_i$ is composed of all points **x** that have $\mathbf{p}_i$ as the nearest sampling point by Euclidean distance $d(\mathbf{x}, \mathbf{p}_i)$, i.e.,

$$\{d(\mathbf{x},\mathbf{p}_i) \leq d(\mathbf{x},\mathbf{p}_j), \forall j \neq i\}, \quad i,j \in \{1,...,b\} \quad (13)$$

The boundaries of Voronoi regions are equidistant from two or more sampling points, while the vertices called Voronoi points are the intersection points ($\mathbf{Q}_i$ in Fig. 2) of the perpendicular bisectors. Triangulation graphs (blue triangles in Fig. 2) are formed by linking each sample point with the points in its neighboring cells. These triangles are characterized by the property that the circumcircle (or circumsphere in higher dimensions) of any triangle (or simplex) does not contain any other points of $\mathbb{P}$ in its interior. This graph is referred to as the Delaunay triangulation, which is the dual graph of the Voronoi diagram, and both structures contain equivalent information in some sense.

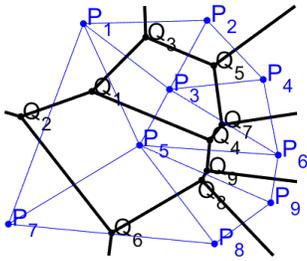

**Fig. 2**. Voronoi diagram

A minimum of $N+1$ initial sample points is required to construct a Voronoi diagram of $N$-dimensional search space. To accurately capture the geometrical characteristics of the search space, $10(N+1)$ initial points are used. The Quickhull algorithm [22] is efficient to compute the convex hull, Voronoi diagram, and Delaunay triangulation. A Delaunay triangulation in $\mathbb{R}^N$ can be derived from a convex hull in $\mathbb{R}^{N+1}$ [22]. The algorithm computes the convex hull in worst-case in $O(b^{(N/2)}/(N/2)!)$ time for $b$ points. Although SDP and SOCP relaxations are polynomial-time solvable in theory, their high computational cost limits scalability for large systems exceeding 500 buses. Furthermore, studies have shown that SDP relaxations are not generally exact, particularly when idealized conditions such as purely resistive networks or small phase angle differences are not satisfied [8]. We use the Qhull general dimension code to generate the Delaunay triangulation, from which the Voronoi diagram is derived, as this approach is computationally more efficient than constructing the Voronoi diagram directly. This Voronoi diagram practically offers several advantages: **(i)** *Systematic search space exploration*, the Voronoi diagram identifies sparsely explored regions by analyzing distances between Voronoi and sample points; **(ii)** *Adaptive refinement and coverage*, each iteration progressively refines the diagram and promotes uniform search space coverage; **(iii)** *Scalability and parallelization*, the Voronoi structure integrates seamlessly with clustering-based partitioning, allowing efficient, parallel computation of local diagrams across clusters while preserving global consistency through boundary sample sharing; and **(iv)** *Feasibility awareness and robustness*, by covering the entire search domain, including disconnected and infeasible regions, the Voronoi diagram effectively captures multiple feasible areas.

### 2.5. Step 4: Identify the tentative optimum and the candidate Voronoi region

The candidate region, in **Step 4**, is defined as the Voronoi region containing the tentative optimal solution, identified as the sample point with the minimum cost function value. Improving the tentative optimum can be approached in different ways. For example, in [20], local function expressions are generated as quadratic polynomials for respective sample points by applying the least squares method, then it is perturbed within triangle comprised by the sample point of Voronoi candidate region and Voronoi points at the border of the region. In cases involving time-intensive cost functions computing, such as in [21], Barycentric interpolation is used to estimate the local function around sample points, then it is refined using linear programming. In the proposed approach, the CPG method with enforced constraint restoration is employed to upgrade the tentative optimum as follows:

As illustrated in Fig. 3, if the tentative optimal solution is identified as $\mathbf{p}_i$, the shaded Voronoi region (polyhedron) is considered the most effective area. To prepare $\mathbf{p}_i$ for integration, its correspondence state variables and functional constraints are first utilized to compute the slack variables using (2), then all variables are used as an initial point of $\psi(\mathbf{x})$. To keep the end point of the integral path within the candidate region, a set of linear inequality constraints defining the candidate region is imposed. Thus, three scenarios arise: 1) the integral curve reaches the boundary of the candidate region, identifying a new point $\mathbf{p}_i^*$ at the boundary as shown in Fig. 3; 2) the curve converges to an asymptotically stable solution, i.e., an optimal solution, within the candidate region; or 3) the initial point $\mathbf{p}_i$ is itself an asymptotically stable solution.

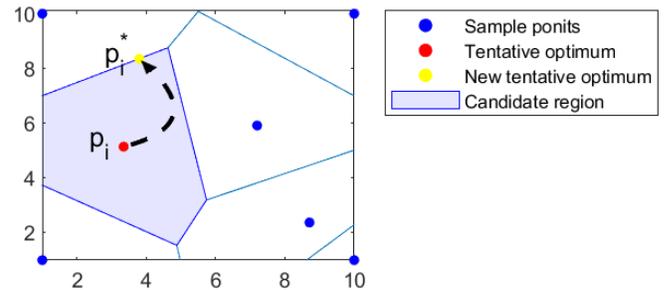

**Fig. 3**. Candidate region and tentative optimum.

The second term of $\psi(\mathbf{x})$ represents the projection of the cost function gradient $\nabla f(\mathbf{x})$ onto the tangent space $\mathbf{T}_X$ at every feasible point, aligning with the PGS system proposed in [3]. This term and the PGS system are identical to $\psi(\mathbf{x})$ on the feasible set, meaning they yield the same solution if initialized at the same point within the feasible set. Consequently, the second term of $\psi(\mathbf{x})$ requires a correction procedure to address any constraint violations when initialized outside the feasible space. The first term of $\psi(\mathbf{x})$ serves as a correction component to restore constraint adherence and is orthogonal to the cost function gradient projection, i.e., the second term of $\psi(\mathbf{x})$, as shown in Fig. 4. If the initial point is infeasible, then correction term typically outweighs the cost function gradient projection term. Consequently, the solution trajectories are first guided toward the feasible space. Once within the feasible set, which is an invariant set of $\psi(\mathbf{x})$, the trajectories progress along it toward an optimal point. Therefore, the enforced constraint restoration feature of $\psi(\mathbf{x})$ allows starting from a point that may not necessarily be feasible as the initial point; however, the integral curve will still converge to an OPF solution.

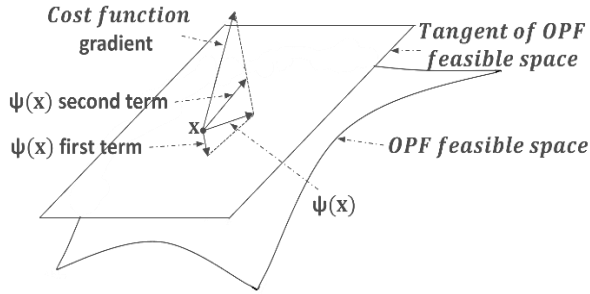

**Fig. 4**. $\psi(\mathbf{x})$ geometry.

The equilibria of $\psi(\mathbf{x})$ correspond exactly to the critical points of the OPF problem, and only the optimal points are asymptotically stable, ensuring that the system naturally converges toward OPF optima. As a result, $\psi(\mathbf{x})$ defines a vector field in which OPF minima act as attractors, guaranteeing convergence to feasible and optimal solutions under the conditions established in Theorems 1 and 2.

*2.6. Step 5: Identify and add new sample points*

In **Step 5**, enhancing the tentative optimum and improving the global approximation fidelity are two primary objectives that can be served by adding new sample points [20]. With each iteration, new sample points are incorporated to refine the Voronoi diagram and better explore the search space. To achieve these objectives, three new sample points are added to the set in each iteration, each possessing distinct characteristics.

The first new sample point is the **integration endpoint,** where the integral path of $\psi(\mathbf{x})$ converges, starting from the tentative optimum as the initial point.

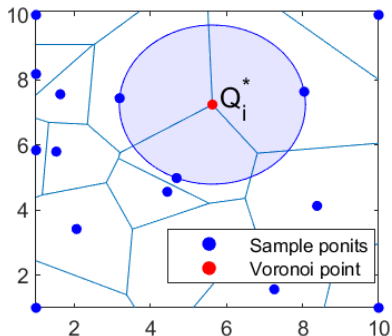

**Fig 5**. Most depopulated region.

The second point is the **farthest Voronoi point** located in the most depopulated region to ensure the global approximation fidelity of the search space. This region, $V_i$, is characterized by a Voronoi point $\mathbf{Q}_i^*$ that has the largest distance $d(\mathbf{Q}_i^*, \mathbf{p}_j)$ to its surrounding sample points $\mathbf{p}_j$, where $i \in \mathbb{L}_i$, $\mathbb{L}_i$ is the set of sample points surrounds $\mathbf{Q}_i \in \mathbb{Q}$, and $\mathbb{Q}$ is the set of Voronoi points. This region is shown in Fig. 5. Notably, that the distances from $\mathbf{Q}_i^*$ to all surrounding sample points are equal, forming a circle, largest empty circle, centered at $\mathbf{Q}_i^*$ with radius of $d(\mathbf{Q}_i^*, \mathbf{p}_j)$.

The third is a **balancer point** placed at the midpoint of the line connecting the first and second points. This can mitigate the uneven populations in regions around the added first and second sample points to prevent a distortion of the Voronoi diagram at the interface between these regions [20].

If the tentative optimum was already used as the initial point for $\psi(\mathbf{x})$ in a previous iteration, the second farthest Voronoi point is selected as the initial point which prevents stalling at local optima and improves coverage of the search space. Finally, new sample points are determined and incorporated into the current sample points set, completing one iteration ensuring the global fidelity of the Voronoi diagram by preserving its fractal structure, while the CPG method accelerates the local search for optimal points. This iterative enrichment enhances global exploration and reduces local trapping, albeit with increased computational cost.

## 3. OPF Formulation

The OPF problem usually minimizes the generation cost by adjusting control variables while satisfying operational constraints. A quadratic cost function $f_i$ is assumed in this paper in terms of each active power generation. Consider a power system network with $N_B$ buses, $N_G$ generators, and $N_L$ lines. The following OPF formulation is considered:

$$\min \quad \sum_{i \in \mathbb{G}} f_i(P_{Gi}) \tag{14}$$

subject to:

$$\begin{cases} V_i \sum_{j \in \mathbb{B}} V_j (G_{ij}\cos(\theta_i - \theta_j) + B_{ij}\sin(\theta_i - \theta_j)) - P_{Gi} + P_{Di} = 0 \\ V_i \sum_{j \in \mathbb{B}} V_j (G_{ij}\sin(\theta_i - \theta_j) - B_{ij}\cos(\theta_i - \theta_j)) - Q_{Gi} + Q_{Di} = 0 \end{cases}, i \in \mathbb{B} \tag{15}$$

$$V_i^{\min} \le V_i \le V_i^{\max}, \quad i \in \mathbb{B} \tag{16}$$

$$P_{Gi}^{\min} \le P_{Gi} \le P_{Gi}^{\max}, \quad Q_{Gi}^{\min} \le Q_{Gi} \le Q_{Gi}^{\max}, \quad i \in \mathbb{G} \tag{17}$$

$$\left|S_l^f\right| \le S_l^{\max}, \quad \left|S_l^t\right| \le S_l^{\max}, \quad l \in \mathbb{L} \tag{18}$$

$$\theta_{slack} = 0 \tag{19}$$

where $\mathbb{B}$, $\mathbb{G}$ and $\mathbb{L}$ are respectively the sets of buses, generator buses and lines. At each bus $i$, $P_{Gi}$ and $Q_{Gi}$ are the real and reactive power injections of the generator, $V_i$ and $\theta_i$ are the voltage magnitude and angle, and $P_{Di}$ and $Q_{Di}$ are the real and reactive power loads, respectively; parameters $G_{ij}$ and $B_{ij}$ are the conductance and susceptance of a transmission line between buses $i$ and $j$; $S_l^f$ and $S_l^t$ are the complex power flows of line $l$ at the from- and to-buses at line $l$, respectively. Note that all variables are continuous. $V_i^{max}$ and $V_i^{min}$ are the upper and lower bounds on variable $V_i$; $P_{Gi}^{max}$ and $P_{Gi}^{min}$ are the upper and lower bounds on variable $P_{Gi}$; $Q_{Gi}^{max}$ and $Q_{Gi}^{min}$ are the upper and lower bounds on variable $Q_{Gi}$; $S_l^{max}$ is the upper bound on variables $S_l^f$ and $S_l^t$. Note that constraints in (18) account for complex power flow differences caused by line losses and impose limits at both ends of the line. Constraint (19) is imposed to set an angle reference. Considering (15) and the functional expressions $S_l$ in (18) are nonlinear equations, the OPF is a nonconvex problem.

The control variables (**u**) include the generators real power injections except the generator at slack bus and the voltage magnitude of the P-V buses and slack bus, thus $\mathbf{u} \in R^{2N_G-1}$, while other variables such as P-Q bus voltage magnitude and angle and generators reactive power injections are the state variables $\mathbf{y} \in R^{2N_B}$ and can be calculated via the power flow equations when the control variables are evaluated. Based on that, functions in (1) have the following dimensions: $f(\mathbf{u},\mathbf{y}): \mathbb{R}^{2N_G+2N_B-1} \to \mathbb{R}$, $\mathbf{e}(\mathbf{u},\mathbf{y}): \mathbb{R}^{2N_G+2N_B-1} \to \mathbb{R}^{2N_B+1}$, and $\mathbf{h}(\mathbf{u},\mathbf{y}): \mathbb{R}^{2N_G+2N_B-1} \to \mathbb{R}^{4N_G+2N_B+2N_L}$. In (3), $n = 6N_G + 4N_B + 2N_L - 1$ and $m = 4N_G + 4N_B + 2N_L + 1$. Thus, the Jacobian matrix $\mathbf{J_g}(\mathbf{x})$ has a full row rank equal to $m$ and the feasible space dimension equals to $2N_G - 1$.

## 4. Case Studies

This section demonstrates how the proposed approach successfully achieves the global OPF solutions for several test systems originally introduced and analyzed in [2], which all have multiple local minima. When working with higher-dimensional data such as those from the IEEE 39-bus and 118-bus systems, partitioning the data into smaller subsets and computing the Voronoi diagram for each subset independently becomes more computationally efficient than computing the Voronoi diagram for all points at once. Clustering algorithms, such as k-Means, can be utilized to partition a dataset into clusters. Once the data is divided into subsets, the Voronoi diagram is computed for each subset independently. To preserve the global structure of Voronoi diagram when identifying the most depopulated area, the distance of Voronoi points that lie outside of a cluster is calculated using the nearest sample points from all clusters, rather than only points from the cluster where they have been constructed. In addition, the Quickhull algorithm systematically partitions data to construct the convex hull and complete Delaunay triangulation, ensuring full coverage of the search space, including sparse regions. While sparse areas yield larger triangles and Voronoi regions, no points are omitted, and accuracy depends on sample distribution rather than algorithmic limitations. However, the nearest sample points from neighboring clusters are incorporated when calculating the distances of Voronoi points, ensuring accuracy and continuity across partitions. Hence, for large-scale systems with thousands of buses, where the search space dimensionality may reach hundreds or thousands, this partitioning enables parallel computation, as each cluster can be processed independently across multiple cores, significantly accelerating the overall computation. For instance, given the worst-case complexity of $O(b^{(N/2)}/(N/2)!)$, dividing the search space into $c$ clusters reduces the overall complexity approximately by $(b^{(N(c-1)/2c)} * (N/2c)!/(N/2)!)$ representing a significant decrease relative to the unpartitioned case. To simplify the construction of the Voronoi diagram, each cluster has a maximum sample points number limit, and all sample points will be re-clustered every iteration. The predefined stopping criterion for global optimality in the current implementation is based on a fixed number of iterations.

We provide a concise summary linking the algorithmic procedures to the entries presented here. The *Tentative Optimum* refers to the sample point with the minimum cost value in each iteration, which serves as the initial point for evaluating $\psi(\mathbf{x})$. The *Second Farthest Point* ($P_i$) is selected as the initial point for $\psi(\mathbf{x})$ when the current tentative optimum has already been used as an initial point in a previous iteration. The *First Added Point* ($S_i$) generated by integrating the initial point in each iteration. Finally, the Balancer points (not explicitly listed) are placed at midpoints between integration endpoints and farthest points to mitigate uneven sampling.

### 4.1. IEEE 9-bus System

This system consists of nine buses, three generators, and nine transmission lines with quadratic cost functions. Therefore, the feasible set is five dimensions: the real output power of the generator at buses 2 and 3, denoted as $P_{G2}$ and $P_{G3}$, along with the voltage magnitudes of buses 1, 2, and 3: $V_1$, $V_2$ and $V_3$. Thus, the Voronoi diagram is five-dimensional. This case has four optimal solutions; the cost of the best local solution is 10% higher than the cost of the global solution. The generators' reactive power outputs lower limit was raised from -300 Mvar to -5 Mvar and they are binding in all optimal solutions.

To test the robustness of the proposed approach, the second best optimal solution or the local optimal solution, $S_{local}$, which occurs at: $f(x) = 3398.03, P_{G2} = 0.648, P_{G3} = 1.178, V_1 = 0.9064, V_2 = 0.9255, V_3 = 0.9326$, along with other 49 random sample points, each with cost functions higher than that of $S_{local}$, were used to constitute a set of 50 initial sample points. Therefore, in the first iteration, the tentative optimum is $S_{local}$. Being an asymptotically stable point, the integral path of $\Psi(x)$ remains at it. To prevent duplicating existing points, the tentative optimum is slightly perturbed within the candidate region serving as the integration endpoint. After three points were added to the initial sample points, $S_{local}$ still has the minimum cost function, i.e., it is the tentative optimum, Fig. 6(a), (b), and (c) illustrate the types of sample points and added sample points of the first, third, and fifth iterations, respectively.

**In iteration 2**, the second farthest point, $P_2$, at ($P_{G2} = 3.0$, $P_{G3} = 1.1914, V_1 = 0.90, V_2 = 1.10, V_3 = 0.90$) is selected as an initial point, the integral curve reach to the first new added point $S_2 : (f(\mathbf{x}) = 4037.2, P_{G2} = 1.733, P_{G3} = 0.102, V_1 = 0.903, V_2 = 0.907, V_3 = 0.913)$.

**In iteration 3**, the second farthest point, $P_3$, at ($P_{G2} = 0.5282$, $P_{G3} = 2.0723, V_1 = 0.920, V_2 = 1.10, V_3 = 0.90$) is selected as an initial point, the integral curve reach to the first new added point $S_3$: ($f(\mathbf{x}) = 5591.8, P_{G2} = 0.3151, P_{G3} = 1.5221, V_1 = 0.9188, V_2 = 0.9862, V_3 = 0.90$).

**In iteration 4**, the second farthest point was again selected as an initial point, $P_4$, which is located at ($P_{G2}$ = 0.0, $P_{G3}$ = 1.9542, $V_1$ = 1.10, $V_2$ = 1.10, $V_3$ = 1.0342). Integrating $P_4$ led to the first added point, $S_4$: ($f(\mathbf{x})$ = 6471.7, $P_{G2}$ = 0.10, $P_{G3}$ = 1.7295, $V_1$ =1.0595, $V_2$ = 1.0687, $V_3$ = 1.0369).

**In iteration 5**, similarly the second farthest point, $P_5$, at ($P_{G2}$ = 0.5470, $P_{G3}$ = 1.1070, $V_1$ = 0.9762, $V_2$ = 1.0550, $V_3$ = 1.0159). The integration endpoint obtained at $S_5$: ($f(\mathbf{x})$ = 3282.9, $P_{G2}$ = 0.8055, $P_{G3}$ = 0.9962, $V_1$ = 0.9064, $V_2$ = 0.9262, $V_3$ = 0.9336), as shown in Fig. 7, making it the tentative optimum in the next iteration since it yields the lowest cost.

**In iteration 6**, integrating $S_5$ led to $S_6$: ($f(\mathbf{x})$ = 3144.6, $P_{G2}$ = 1.0316, $P_{G3}$ = 0.7904, $V_1$ = 0.9080, $V_2$ = 0.9238, $V_3$ = 0.9338)

**Finally, in iteration 7**, utilizing $S_6$ as tentative optimum, the global optimal solution is obtained at ($f(\mathbf{x})$ = 3087.8, $P_{G2}$ = 1.2537, $P_{G3}$ = 0.5703, $V_1$ = 0.9095, $V_2$ = 0.9218, $V_3$ = 0.9388), as illustrated by Fig. 7.

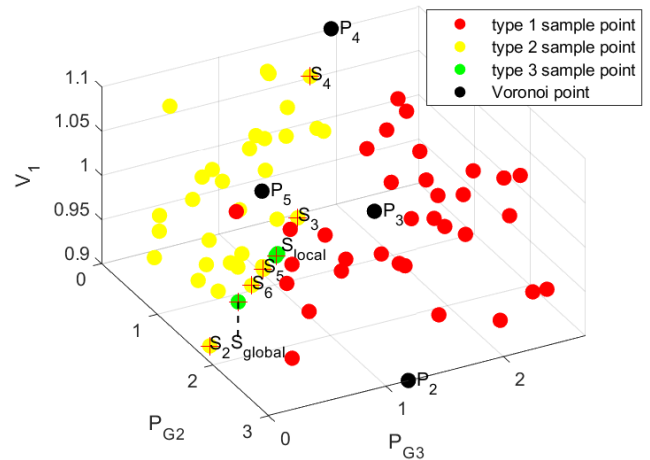

**Fig. 7.** Case9mod: global OPF is achieved at iteration 7.

### 4.2. IEEE 39-bus System

The IEEE 39-bus system ("case39mod2") has 39 buses, 46 transmission lines and 10 generators. Consequently, the feasible set is 19-dimensional, comprising nine real output power variables and ten bus voltage magnitudes. This case features 16 optimal solutions with very different generation levels; however, they are all within 0.5% in objective value when only linear cost coefficients are used. Additionally, it has disconnected feasible region. 190 initial sample points are utilized to construct the Voronoi diagram, the points are partitioned into clusters each containing a maximum of 25 points. The best local optimum, $S_{local}$, as detailed in Table 1, is included in the initial sample points set.

In **iteration 1**, since $S_{local}$ has the minimum cost, it is selected as the tentative optimum, $S_{local}$. As an asymptotically stable equilibrium point, the integral curve remained unchanged. After adding three points to the initial sample set, $S_{local}$ continued to exhibit the lowest cost function. As a result, during **iteration 2**, the second farthest point, $P_2$, was chosen as the initial point, leading the integral curve to converge to the first newly added point in that iteration, $S_2$. Similarly, in **iteration 3**, the second farthest point, $P_3$, was selected as the initial point. Integrating $P_3$ produced $S_3$, the first point added during the third iteration. In **iteration 4**, $P_4$, the second farthest point was used as the initial point. Integration of $P_4$ resulted in $S_4$, which became the tentative optimum for the subsequent iteration. From **iterations 5 to 7**, all tentative optimum points consistently yielded the lowest cost functions. Finally, the global optimum, $S_{global}$, was reached after seven iterations. The repeated selection of the second-farthest point over three consecutive iterations indicates that the tentative optimum has stalled, prompting the algorithm to redirect its search toward sparsely sampled regions of the solution space. The integration endpoints ($S_2$, $S_3$) illustrate how the method expands exploration into alternative areas of the feasible region. In contrast, the sequence ($S_4$, $S_5$, $S_6$, $S_{global}$) demonstrates a monotonic progression toward the global optimum.

### 4.3. IEEE 118-bus System

This system ("case118mod") is composed of 118 buses, 54 generators, and 186 transmission lines which results in a 107-

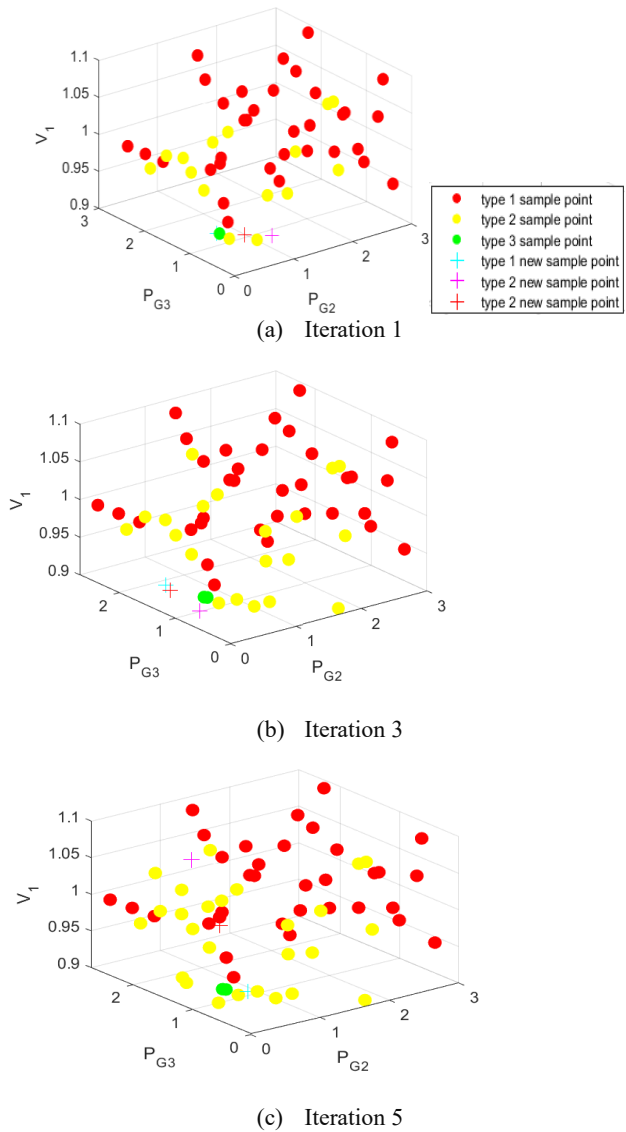

**Fig. 6.** Case9mod: sample points added in iterations 1, 3 and 5.

dimensional feasible set consisting of 53 real output power and 54 voltage magnitude. Three optimal solutions where the cost of the local solution is 38% greater than the global solution. This case exhibits disconnected feasible regions, too. 1070 initial sample points are utilized to construct the Voronoi diagram, the points are partitioned into clusters each containing a maximum of 110 points. The best local optimum, $S_{local}$, as detailed in Table 2, is included in the initial sample points set. Table 2 maps the algorithm steps to key sample points generated during the iterations. The progression of cost values across these points illustrates the method's systematic convergence from local to global solutions.

During **iteration 1**, $S_{local}$, was designated as the tentative optimum since it has the lowest cost. As an asymptotically stable equilibrium, the integral curve remained unchanged. After incorporating three additional points into the initial sample set, $S_{local}$ still maintained the minimum cost function. Consequently, in **iteration 2**, the second farthest point, $P_2$, was selected as the initial point, guiding the integral curve toward the first newly introduced point in that iteration, $S_2$. The process was repeated in **iteration 3**, where $P_3$ was chosen as the initial point, leading to the integration result $S_3$, the first point added in this iteration. In **iteration 4,** the second farthest point, $P_4$, was selected as the initial point. Integrating $P_4$ produced $S_4$. In **iteration 5**, $P_5$, the second farthest point was used as the initial point, and its integration produced $S_5$, which then became the tentative optimum for the next iteration. Throughout **iterations 6 to 8**, every tentative optimum consistently achieved the lowest cost function. Ultimately, after eight iterations, the global optimum, $S_{global}$, was successfully attained.

### 4.4. Sensitivity analysis

To evaluate how the proposed approach responds to variations in sample point distribution, a new 9-bus system case study is designed. In this setup, the initial sample points are placed only in the upper half of the search space by raising the lower bounds of the control variables to their midpoints. Figure 8 illustrates both the initial and newly added points, showing that the global optimal solution is obtained after six iterations. During iterations 2 to 5, the second farthest point is selected as initial integration point.

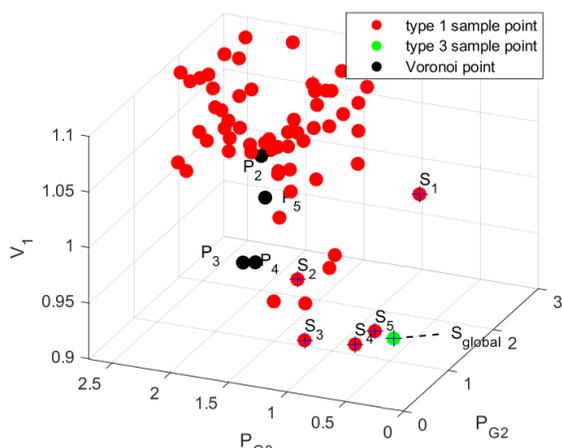

**Fig. 8.** Case9mod: global OPF is achieved at iteration 6

The sensitivity to initial sample points distribution is mitigated by the iterative enrichment process. This strategy progressively balances the distribution of points and ensures coverage of the search space through Voronoi regions around optimal points. In addition, when sample points are unevenly distributed or certain areas are sparsely populated, the corresponding Voronoi region becomes larger, allowing the CPG integration path to extend further and terminate at the region's boundary.

To evaluate the impact of the penalty parameter $c$ in (12) for category II points, various values were tested in the 9-bus system, using the same initial sample points as described in Section 4.1. The parameter $c$ was set to 500, 1000, 3000, 5000, and 7000. The three highest values produced identical added points and converged in the same number of iterations. On the other hand, when ccc is less than 1000, the tentative optimum selected in the second iteration is a type II point with a lower cost than the local optimal solution, $S_{local}$, due to the small penalty term. However, this point is not asymptotically stable, and the integration path eventually leaves it, leading to convergence to the global optimum after seven iterations.

### 4.5. Comparison of the proposed method with IPOPT

This subsection summarizes the comparisons of the proposed Voronoi diagram-based OPF method with the IPOPT solver on four test systems. Including the three systems that have been studied above and a new case, "WB5mod", that is studied in detail below. Table 3 provides a consolidated comparison of the objective value, computation time, and optimality gap. As shown in the table, the proposed method achieves an optimality gap below 0.001% in all cases, whereas IPOPT failed to find the global optimum for the "WB5mod" case starting from a flat initialization. In terms of runtime, the proposed method can be slower, especially for larger systems. This is mainly due to the time required to construct and compute the Voronoi diagrams.

In the following, we detail the implementation for the "WB5mod" case [15]. It has 5 buses, 6 transmission lines and 2 generators. The IPOPT solver converges to local optimal solution $S_{local}$: ($f(\mathbf{x}) = 161921.2$, $P_{G2} = 0.3982$, $V_1 = 1.10$, $V_5 = 0.9232$, from flat start. The feasible set is 3-dimensional, comprising one real output power variables and two bus voltage magnitudes. 30 initial sample points are utilized to start the proposed approach, $S_{local}$ is amongst them. Therefore, in the first iteration, the tentative optimum is $S_{local}$. Being an asymptotically stable point, the integral path of $\Psi(x)$ remains at it. To prevent duplicating existing points, the tentative optimum is slightly perturbed within the candidate region serving as the integration endpoint. In the second iteration, the second farthest point, $P_2$, at ($P_{G2} = 4.225$, $V_1 = 0.9080$, $V_5 = 0.9408$) is selected as an initial point, the integral curve reach to the first new added point $S_2$: ($f(\mathbf{x}) = 191995.3$, $P_{G2} = 0.5659$, $V_1 = 1.0209$, $V_5 = 0.9205$). In iteration 3, the tentative optimum is $S_{local}$ again, the second farthest point, $P_3$, at ($P_{G2} = 1.0953$, $V_1 = 0.9071$, $V_5 = 1.6047$) is selected as an initial point, the integral curve reach to the first new added point $S_3$: ($f(\mathbf{x}) = 147251.9$, $P_{G2} = -.110$, $V_1 = 1.1086$, $V_5 = 0.9588$). $S_3$ is selected as the tentative optimum for the first time, hence it serves as initial integral point, leading to point $S_4$: ($f(\mathbf{x})$ = 139873.3, $P_{G2}$ = 0.0, $V_1$ = 1.0892, $V_5$ = 0.9385) making it the tentative optimum in the iteration 5. Where the integration path did not leave $S_4$ indicating optimal solution. In iteration 6, the second farthest point, $P_6$: ($P_{G2}$ = 0.4725, $V_1$ = 1.0331, $V_5$ = 1.0747) is selected as an initial point, the integral curve reach to $S_6$: ($f(\mathbf{x})$ = 562799.4, $P_{G2}$ = 0.7440, $V_1$ = 0.9991, $V_5$ = 1.0726). In final iteration, 7, the integral path reached $S_7$: ($f(\mathbf{x})$ = 276225.1, $P_{G2}$ = 0.8336, $V_1$ = 0.90, $V_5$ = 0.9173). from $P_7$: ($P_{G2}$ = 3.6154, $V_1$ = 1.0646, $V_5$ = 0.9856).

## 5. Conclusion

This paper presents a novel approach for the global OPF solution by constructing the Voronoi diagram about the geometric structure of the search space. Also, the CPG method is utilized to establish the connection between the asymptotically stable equilibrium points of an equivalent autonomous system and the local optima of the OPF, leveraging its constraints restoration capability to initialize the autonomous system not at necessarily a feasible point. To iteratively refine the Voronoi diagram, the proposed approach adds three new sample points at each iteration including the enhanced tentative optimal point determined by the CPG, the farthest Voronoi point, and the balancer point. A numerical procedure is also proposed to efficiently compute the global optimal solution for OPF problems. The effectiveness of the proposed approach for global OPF solutions has been validated using IEEE test cases, successfully achieving global optimality in benchmark systems exhibiting multiple optima presented in [2]. In contrast, semidefinite relaxations of the OPF problem have shown limited success in closing the duality gap, as seven of the eight test cases containing local optima yielded nonzero duality gaps under certain parameter settings [2]. Future work will focus on integrating derivative-free local solvers to handle non-differentiable objective functions and constraints, as well as extending the framework to include discrete control variables. In addition, efforts will be directed toward developing a more efficient partitioning strategy to enable real-time AC OPF computations. While methods such as moment-SOS relaxation [23], which has computational requirements substantially higher than those of SDP relaxations [8], and branch-and-bound [24] demonstrate robustness and theoretical optimality guarantees, their scalability to large-scale systems remains limited compared to the proposed Voronoi-based approach.

**Table 1** Key points of "case39mod2" power injections & bus voltages

|  | $S_{local}$ | $P_2$ | $S_2$ | $P_3$ | $S_3$ | $P_4$ | $S_4$ | $S_5$ | $S_6$ | $S_{global}$ |
|---|---|---|---|---|---|---|---|---|---|---|
| $P_{G30}$ | 0.00 | 10.4000 | 10.4000 | 10.4000 | 10.4000 | 0.4609 | 1.8625 | 5.7967 | 7.2942 | 7.1179 |
| $P_{G32}$ | 725.00 | 0.0000 | 0.0001 | 0.0000 | 0.0000 | 6.8597 | 7.2320 | 6.6948 | 1.1329 | 0.0004 |
| $P_{G33}$ | 202.87 | 6.5200 | 6.5200 | 6.5193 | 2.6807 | 2.9348 | 2.0120 | 1.7635 | 1.8597 | 1.8893 |
| $P_{G34}$ | 226.90 | 0.3066 | 2.0899 | 5.0798 | 2.4393 | 2.4802 | 2.1333 | 2.1700 | 2.2393 | 2.2501 |
| $P_{G35}$ | 399.11 | 0.0000 | 0.0000 | 0.0010 | 4.8816 | 4.4829 | 4.2993 | 4.0969 | 4.7713 | 4.6585 |
| $P_{G36}$ | 126.15 | 0.0000 | 0.0001 | 0.0000 | 0.0000 | 1.2040 | 1.0245 | 0.7183 | 0.3495 | 0.5437 |
| $P_{G37}$ | 454.24 | 5.6400 | 4.5977 | 2.0071 | 1.7978 | 2.5855 | 2.2597 | 2.3904 | 2.3940 | 2.4033 |
| $P_{G38}$ | 59.320 | 0.0000 | 0.0252 | 0.0001 | 2.8167 | 3.2573 | 2.6863 | 2.4641 | 2.5641 | 2.5783 |
| $P_{G39}$ | 949.53 | 0.0000 | 1.4705 | 0.0000 | 0.0009 | 3.2286 | 2.4563 | 2.8778 | 3.3560 | 3.4820 |
| $V_{30}$ | 1.048 | 0.9500 | 1.0359 | 0.9974 | 1.0266 | 1.0306 | 1.0481 | 1.0437 | 1.0409 | 1.0412 |
| $V_{31}$ | 0.950 | 0.9500 | 0.9523 | 1.0500 | 0.9816 | 0.9736 | 0.9500 | 0.9500 | 0.9500 | 0.9500 |
| $V_{32}$ | 0.974 | 0.9500 | 1.0060 | 1.0500 | 1.0111 | 0.9845 | 0.9664 | 0.9768 | 0.9974 | 0.9957 |
| $V_{33}$ | 0.950 | 1.0446 | 0.9844 | 0.9570 | 0.9889 | 1.0000 | 0.9500 | 0.9500 | 0.9500 | 0.9500 |
| $V_{34}$ | 0.950 | 1.0500 | 1.0159 | 0.9500 | 0.9862 | 0.9873 | 0.9500 | 0.9500 | 0.9500 | 0.9500 |
| $V_{35}$ | 0.950 | 1.0131 | 0.9921 | 0.9501 | 0.9904 | 0.9686 | 0.9500 | 0.9500 | 0.9500 | 0.9500 |
| $V_{36}$ | 0.994 | 0.9500 | 1.0281 | 1.0500 | 1.0472 | 0.9766 | 0.9927 | 0.9943 | 0.9952 | 0.9951 |
| $V_{37}$ | 1.022 | 0.9500 | 1.0245 | 0.9500 | 1.0162 | 1.0263 | 1.0244 | 1.0204 | 1.0192 | 1.0193 |
| $V_{38}$ | 0.961 | 0.9500 | 0.9913 | 0.9500 | 0.9847 | 0.9912 | 0.9692 | 0.9670 | 0.9668 | 0.9667 |
| $V_{39}$ | 1.025 | 0.9500 | 1.0359 | 0.9932 | 1.0005 | 1.0045 | 1.0052 | 1.0128 | 1.0195 | 1.0191 |
| Cost | 942.95 | - | 946.91 | - | 944.31 | - | 942.66 | 942.11 | 941.79 | 941.738 |

Table 2 Key Points of "case118mod" Power Injections

| | $S_{local}$ | $P_2$ | $S_2$ | $P_3$ | $S_3$ | $P_4$ | $S_4$ | $P_5$ | $S_5$ | $S_6$ | $S_7$ | $S_{global}$ |
|---|---|---|---|---|---|---|---|---|---|---|---|---|
| $P_{G1}$ | 0.0000 | 0.0000 | 0.0000 | 0.0000 | 0.0000 | 5.4613 | 3.6511 | 0.0000 | 0.0000 | 0.0000 | 0.0038 | 0.2555 |
| $P_{G4}$ | 0.0000 | 0.0000 | 0.0000 | 7.0000 | 2.9473 | 7.0000 | 3.4148 | 2.2323 | 0.7318 | 0.5888 | 0.4386 | 0.0000 |
| $P_{G6}$ | 0.0000 | 7.0000 | 3.0118 | 0.0000 | 0.0000 | 0.0000 | 0.0000 | 0.0000 | 0.0000 | 0.0000 | 0.0000 | 0.0000 |
| $P_{G8}$ | 0.0000 | 7.0000 | 5.5416 | 4.5068 | 4.4383 | 3.6513 | 3.2218 | 2.2298 | 0.2026 | 0.1897 | 0.1559 | 0.0000 |
| $P_{G10}$ | 1.2360 | 0.0000 | 0.0000 | 12.5892 | 4.8879 | 14.7575 | 6.1240 | 4.4171 | 4.6514 | 4.1237 | 4.1186 | 4.0155 |
| $P_{G12}$ | 0.2645 | 12.9500 | 1.6737 | 0.0000 | 0.0000 | 0.5151 | 0.0005 | 0.0001 | 0.0099 | 0.7713 | 0.8856 | 0.8572 |
| $P_{G15}$ | 0.0000 | 7.0000 | 2.1359 | 0.0000 | 0.0000 | 1.2652 | 2.3442 | 0.0000 | 0.3293 | 0.3662 | 0.3850 | 0.1997 |
| $P_{G18}$ | 0.0000 | 7.0000 | 0.0000 | 7.0000 | 2.3412 | 4.2633 | 2.1442 | 0.1515 | 0.0011 | 0.0016 | 0.0021 | 0.1233 |
| $P_{G19}$ | 0.0000 | 7.0000 | 1.9682 | 0.0000 | 0.0000 | 0.0000 | 0.0000 | 0.0000 | 0.0000 | 0.0000 | 0.0000 | 0.2095 |
| $P_{G24}$ | 0.0000 | 7.0000 | 2.3801 | 0.0000 | 0.0000 | 0.0000 | 0.0001 | 0.0000 | 1.4025 | 0.5664 | 0.3495 | 0.0000 |
| $P_{G25}$ | 0.0000 | 1.1052 | 1.9863 | 22.4000 | 3.3591 | 22.4000 | 3.0049 | 6.3789 | 2.0931 | 2.0153 | 1.9849 | 1.9434 |
| $P_{G26}$ | 0.0000 | 3.5375 | 3.8094 | 0.0000 | 0.0000 | 9.4151 | 4.1219 | 0.0002 | 2.8990 | 2.8697 | 2.8601 | 2.8036 |
| $P_{G27}$ | 0.0000 | 0.0000 | 0.0000 | 7.0000 | 3.2686 | 7.0000 | 0.0004 | 1.9851 | 1.1378 | 0.7934 | 0.5853 | 0.1084 |
| $P_{G31}$ | 0.0000 | 0.0000 | 0.0000 | 0.0000 | 0.0000 | 7.4900 | 0.1273 | 0.0001 | 0.0645 | 0.0735 | 0.0727 | 0.0725 |
| $P_{G32}$ | 0.0000 | 7.0000 | 2.7198 | 7.0000 | 3.2471 | 6.1615 | 0.0002 | 7.0000 | 1.2106 | 0.8552 | 0.6434 | 0.1542 |
| $P_{G34}$ | 0.0000 | 7.0000 | 4.7948 | 3.4911 | 3.2993 | 0.9250 | 3.9916 | 0.7720 | 0.0608 | 0.0585 | 0.0696 | 0.0362 |
| $P_{G36}$ | 0.0000 | 7.0000 | 1.9838 | 0.0000 | 0.0000 | 0.0000 | 0.0000 | 0.0000 | 0.1896 | 0.1733 | 0.1945 | 0.0954 |
| $P_{G40}$ | 0.0000 | 0.0000 | 0.0000 | 7.0000 | 2.2740 | 4.2612 | 0.0944 | 7.0000 | 0.1376 | 0.2408 | 0.4120 | 0.4880 |
| $P_{G42}$ | 0.0000 | 7.0000 | 3.8779 | 7.0000 | 3.3959 | 7.0000 | 0.0045 | 7.0000 | 0.0020 | 0.0047 | 0.0121 | 0.4107 |
| $P_{G46}$ | 0.1754 | -0.0000 | 0.0000 | 8.3300 | 0.7892 | 0.2105 | 0.0569 | 0.1697 | 0.0059 | 0.1721 | 0.2048 | 0.1906 |
| $P_{G49}$ | 1.7159 | -0.0000 | 0.0000 | 0.0000 | 0.0004 | 21.2800 | 0.2715 | 0.0001 | 0.0925 | 1.6645 | 2.0546 | 1.9374 |
| $P_{G54}$ | 0.4460 | 10.3600 | 0.5447 | 10.3600 | 0.9933 | 0.6140 | 0.0000 | 0.4449 | 0.0538 | 0.4612 | 0.5377 | 0.4952 |
| $P_{G55}$ | 0.0000 | 1.0001 | 1.0139 | 0.0000 | 0.0000 | 0.0000 | 0.2010 | 0.0000 | 0.0137 | 0.0142 | 0.0246 | 0.3158 |
| $P_{G56}$ | 0.0000 | 0.0000 | 0.0000 | 7.0000 | 5.4008 | 1.7291 | 0.0000 | 2.4682 | 0.0146 | 0.0149 | 0.0253 | 0.3201 |
| $P_{G59}$ | 1.2690 | 0.0000 | 0.0000 | 0.2801 | 1.7772 | 0.0000 | 0.0744 | 0.0001 | 0.0806 | 1.3007 | 1.6056 | 1.4967 |
| $P_{G61}$ | 1.2310 | 0.0000 | 0.0000 | 18.2000 | 2.1187 | 0.9758 | 0.0422 | 1.2210 | 0.0000 | 0.0007 | 0.2948 | 1.4850 |
| $P_{G62}$ | 0.0000 | 7.0000 | 0.7775 | 0.0000 | 0.0000 | 0.0000 | 0.0067 | 0.0254 | 0.0005 | 0.0003 | 0.0002 | 0.0000 |
| $P_{G65}$ | 2.8676 | 5.7107 | 1.3956 | 34.3700 | 2.6307 | 0.1024 | 0.0086 | 0.8136 | 0.0099 | 1.3443 | 3.6951 | 3.5292 |
| $P_{G66}$ | 2.8811 | 34.4400 | 2.0332 | 5.4913 | 5.7831 | 0.0000 | 0.0709 | 0.0000 | 0.0485 | 2.3107 | 3.6540 | 3.4978 |
| $P_{G70}$ | 4.2659 | 0.0000 | 0.0000 | 7.0000 | 0.0008 | 0.0021 | 0.0025 | 0.0006 | 0.0000 | 0.0000 | 0.0000 | 0.0000 |
| $P_{G72}$ | 7.0000 | 0.0000 | 0.0000 | 4.6752 | 1.6378 | 0.0004 | 0.0000 | 0.0003 | 6.9908 | 2.6146 | 0.7290 | 0.0000 |
| $P_{G73}$ | 2.9085 | 3.2187 | 0.6658 | 0.0000 | 0.0000 | 17.0051 | 0.0282 | 0.0000 | 0.0000 | 0.0000 | 0.0000 | 0.0000 |
| $P_{G74}$ | 2.1951 | 2.9973 | 0.2758 | 3.2699 | 0.0013 | 0.0000 | 0.0011 | 0.0000 | 0.0002 | 0.0001 | 0.0001 | 0.1699 |
| $P_{G76}$ | 1.3833 | 6.0728 | 4.5385 | 0.0000 | 0.0000 | 2.1526 | 0.3226 | 0.0000 | 0.0019 | 0.0010 | 0.0018 | 0.2229 |
| $P_{G77}$ | 0.7431 | -0.0000 | 0.0000 | 7.0000 | 0.8812 | 0.2179 | 0.0226 | 0.0625 | 6.9710 | 3.8608 | 1.2836 | 0.0000 |
| $P_{G80}$ | 4.8576 | 4.8555 | 3.0997 | 40.3900 | 3.3090 | 4.0914 | 0.5414 | 3.1537 | 3.2169 | 4.2370 | 4.5835 | 4.3133 |
| $P_{G85}$ | 0.0000 | 7.0000 | 2.4028 | 0.0000 | 0.0000 | 0.0000 | 2.4813 | 0.0000 | 0.4279 | 0.1955 | 0.1109 | 0.0000 |
| $P_{G87}$ | 0.0392 | 7.2800 | 0.0831 | 7.2800 | 0.0927 | 4.4785 | 0.0726 | 0.0257 | 0.0428 | 0.0387 | 0.0381 | 0.0363 |
| $P_{G89}$ | 5.3733 | 7.3645 | 6.7842 | 5.5637 | 7.9514 | 2.8764 | 7.8131 | 3.0542 | 5.6767 | 5.1589 | 5.2115 | 5.0290 |
| $P_{G90}$ | 0.0000 | 2.1545 | 1.8457 | 0.0000 | 0.0000 | 0.0000 | 2.7497 | 0.0000 | 0.0000 | 0.0000 | 0.0000 | 0.0000 |
| $P_{G91}$ | 0.0000 | 3.5414 | 1.8118 | 7.0000 | 2.8007 | 7.0000 | 0.0000 | 0.0000 | 1.0228 | 0.4962 | 0.2782 | 0.0000 |
| $P_{G92}$ | 0.0000 | -0.0000 | 0.0000 | 0.0000 | 0.0000 | 0.0000 | 0.0000 | 5.1023 | 0.0000 | 0.0000 | 0.0000 | 0.0000 |
| $P_{G99}$ | 0.0397 | 7.0000 | 2.0792 | 7.0000 | 2.3336 | 7.0000 | 0.0000 | 1.2546 | 0.3392 | 0.1748 | 0.1157 | 0.0000 |
| $P_{G100}$ | 2.4509 | 0.0000 | -0.0000 | 5.9223 | 4.6210 | 0.9356 | 3.5233 | 7.1652 | 2.3514 | 2.2652 | 2.3170 | 2.3143 |
| $P_{G103}$ | 0.3967 | -0.0000 | 0.0000 | 0.0000 | 0.0000 | 0.0000 | 0.6228 | 0.0000 | 0.4244 | 0.3661 | 0.3732 | 0.3826 |
| $P_{104G}$ | 0.1742 | 7.0000 | 2.3815 | 0.0000 | 0.0000 | 0.0000 | 5.0278 | 0.0000 | 1.7813 | 1.1666 | 0.6023 | 0.0000 |
| $P_{G105}$ | 0.2475 | 0.0000 | 2.4358 | 7.0000 | 4.4880 | 7.0000 | 0.0000 | 7.0000 | 2.1463 | 2.0000 | 1.9999 | 0.0535 |
| $P_{G107}$ | 0.4011 | 7.0000 | 2.2494 | 0.0000 | 0.0000 | 0.0000 | 0.0001 | 0.0000 | 0.0000 | 0.0000 | 0.0000 | 0.2911 |
| $P_{G110}$ | 0.1912 | 0.0000 | 0.0000 | 0.0000 | 0.0000 | 0.0000 | 0.0002 | 0.0000 | 0.0001 | 0.0002 | 0.0002 | 0.0714 |
| $P_{G111}$ | 0.3564 | -0.0000 | 0.0000 | 0.0000 | 0.0000 | 9.5200 | 0.0018 | 0.0000 | 0.0007 | 0.2361 | 0.3597 | 0.3525 |
| $P_{G112}$ | 0.4304 | 7.0000 | 2.5577 | 0.0000 | 0.0000 | 0.0000 | 0.0007 | 0.0000 | 0.0003 | 0.0011 | 0.0023 | 0.3652 |
| $P_{G113}$ | 0.0000 | 0.0000 | 0.0000 | 0.0000 | 0.0000 | 0.0000 | 0.0000 | 0.0000 | 0.7009 | 0.5501 | 0.4296 | 0.0000 |
| $P_{G116}$ | 0.0000 | 7.0000 | 0.0501 | 7.0000 | 0.0281 | 0.0000 | 0.0000 | 1.5154 | 0.0000 | 0.0000 | 0.0000 | 0.0000 |
| Cost | 177984 | - | 329599 | - | 398727 | | 212184 | - | 177715 | 146161 | 138376 | 129625 |

Table 3 Comparison between Voronoi-OPF and IPOPT

|  | Voronoi-OPF | | | IPOPT | | |
|---|---|---|---|---|---|---|
|  | Cost ($/hr) | Time (sec) | Optimality gap (%) | Cost ($/hr) | Time (sec) | Optimality gap (%) |
| WB5mod | 139873.3 | 0.01125 | < 0.001 | 161921.1 | 0.0209 | 15.76 |
| Case9mod | 3087.8 | 0.0919 | < 0.001 | 3087.8 | 0.0780 | < 0.001 |
| Case39mod2 | 941.7 | 0.3325 | < 0.001 | 941.7 | 0.1024 | < 0.001 |
| Case118mod | 129625 | 10.9387 | < 0.001 | 129625 | 0.1301 | < 0.001 |